\documentclass{amsart}
\usepackage{amsmath, amssymb}
\usepackage[all]{xy}

\usepackage{hyperref}%[colorlinks=true]

\theoremstyle{plain}
\newtheorem*{thm*}{Theorem}
\newtheorem*{rem*}{Remark}

\newtheorem{thm}{Theorem}[section]
\newtheorem{cor}[thm]{Corollary}
\newtheorem{defi}[thm]{Definition}
\newtheorem{prop}[thm]{Proposition}
\newtheorem{lm}[thm]{Lemma}

\newtheorem{claim*}{Claim}

\newtheorem{rem}[thm]{Remark}
\newtheorem{nota}[thm]{Notation}
\numberwithin{equation}{thm}

\newcommand{\F}{{\mathcal{F}}}
\newcommand{\E}{{\mathcal{E}}}
\newcommand{\FF}{{\mathbb{F}_2}}

\newcommand{\Eq}{{\ensuremath{\mathcal{E}_{q}}}}
\newcommand{\Eqd}{{\mathcal{E}_{q}^{\mathrm{deg}}}}
\newcommand{\Tq}{{\ensuremath{\mathcal{T}_{q}}}}

\newcommand{\Fq}{{\ensuremath{\mathcal{F}_{quad}}}}
\newcommand{\Gq}{{\ensuremath{\mathcal{F}_{iso}}}}

\newcommand{\Mn}{{\ensuremath{\mathrm{Mix}_{V, D,\eta}}}}
\newcommand{\M}[1]{{\ensuremath{\mathrm{Mix}_{#1}}}}
\newcommand{\Mab}{{\ensuremath{\mathrm{Mix}_{\alpha,\beta}}}}
\newcommand{\m}[1]{{\ensuremath{\mathrm{\Sigma}_{#1}}}}
\newcommand{\mab}{{\ensuremath{\mathrm{\Sigma}_{\alpha,\beta}}}}

\begin{document}

\title{The mixed functors of the category $\mathcal{F}_{quad}$:
a first study}

\author{Christine Vespa}

\address{Ecole Polytechnique F\'ed\'erale de Lausanne, Institut de G\'eom\'etrie, Alg\`ebre et Topologie, Lausanne, Switzerland.}
\email{christine.vespa@epfl.ch}

\date{\today}

\begin{abstract}
In previous work, we defined the category of functors $\Fq$, associated to $\FF$-vector spaces equipped with a nondegenerate quadratic form.  
In this paper, we define a special family of objects in the category $\Fq$, named the mixed functors. We give the complete decompositions of two elements of this family that give rise to two new infinite families of simple objects in the category $\Fq$. 

$\  $

$\  $
\end{abstract}

\maketitle

%\tableofcontents
In $1993$, Henn, Lannes and Schwartz established in \cite{HLS1} a very strong relation between the Steenrod algebra and the category $\F(p)$ of functors from the category $\E^f$ of finite dimensional $\mathbb{F}_p$-vector spaces to the category $\E$ of
all $\mathbb{F}_p$-vector spaces, where $\mathbb{F}_p$ is the prime
field with $p$ elements. To be more precise, they study the category $\mathcal{U}$ of unstable modules over the Steenrod algebra localizated away from the nilpotent unstable modules $\mathcal{N}il$; they exhibit an equivalence between the quotient category $\mathcal{U}/\mathcal{N}il$ and a full subcategory of the category of functors $\mathcal{F}(p)$. This equivalence is very useful and allows  several important topological results to be derived from algebraic results in the category $\F(p)$. For a recent interesting application of this equivalence to the cohomology of Eilenberg MacLane spaces, we refer the reader to the results obtained by Powell in \cite{PowellMZ}.

An important algebraic motivation for the particular interest in the category $\F(p)$ follows from the link with the modular representation theory and the cohomology  of finite general linear groups. Namely, the evaluation of a functor $F$, object in $\F(p)$, on a finite dimensional vector space $V$  is a $\mathbb{F}_p[GL(V)]$-module. A fundamental result obtained by Suslin in the appendix of \cite{FFSS} and, independently, by Betley in \cite{Betley} relates the calculation of extension groups in the category $\F(p)$ with certain stable cohomology groups of general linear groups. 

It is natural to seek to construct other categories of functors that play a similar r\^ole for other families of algebraic groups and, in particular, for the orthogonal groups.

In \cite{math.AT/0606484}, we constructed the functor category $\Fq$, which has some good properties as a candidate for the orthogonal group over the field with two elements. For instance, the evaluation functors give rise to a coefficient system that allows us to define a system of homology groups. We obtained, in \cite{math.AT/0606484}, two families of simple objects in $\Fq$ related, respectively, to general linear groups and to orthogonal groups. The purpose of this paper is to define a new family of objects in the category $\Fq$, named the mixed functors, which give rise to new simple objects of $\Fq$. The mixed functors are subfunctors of a tensor product between a functor coming from the category $\F:= \F(2)$ and a functor coming from the subcategory $\Gq$ of $\Fq$ defined in \cite{math.AT/0606484}. The structure of the mixed functors is very complex, hence it is difficult to give explicit decompositions in general. However, we give the complete decompositions of two significant elements of this family: the functors $\M{0,1}$ and $\M{1,1}$. These two mixed functors play a central r\^{o}le in the forthcoming paper \cite{Vespa3} concerning the decompositions of the standard projective objects $P_{H_0}$ and $P_{H_1}$ of $\Fq$. We prove in  \cite{Vespa3} that these mixed functors are direct summands of $P_{H_0}$ and $P_{H_1}$. The decomposition of $\M{0,1}$ and $\M{1,1}$ represents a further step in our project to classify the simple objects of this category. 

Recall that in \cite{math.AT/0606484}  we constructed two families of simple objects in
$\Fq$. The first one is obtained using the fully faithful, exact functor $\iota: \F \rightarrow
\Fq$, which preserves simple objects. By \cite{KuhnII}, the simple
objects in $\F$ are in one-to-one correspondence with the irreducible
representations of general linear groups. The second family is obtained
using the fully-faithful, exact functor $\kappa: \Gq \rightarrow
\Fq$ which preserves simple objects, where $\Gq$ is
equivalent to the product of the categories of modules over the
orthogonal groups. The results of this paper are summarized in the following theorem.
\begin{thm*}
Let $\alpha$ be an element in $\{0,1 \}$.
\begin{enumerate}
\item
The functor $\M{\alpha,1}$ is infinite.
\item
There exists a subfunctor $\m{\alpha,1}$ of $\M{\alpha,1}$ such that we have the following  short exact sequence
$$0 \rightarrow \m{\alpha,1} \rightarrow \M{\alpha,1} \rightarrow \m{\alpha,1} \rightarrow 0.$$
\item
The functor $\m{\alpha,1}$ is uniserial with unique composition series given by the decreasing filtration given by the subfunctors $k_d \m{\alpha,1}$ of $\m{\alpha,1}$:
$$\ldots \subset k_d \m{\alpha,1} \subset \ldots \subset k_1
\m{\alpha,1} \subset k_0 \m{\alpha,1}= \m{\alpha,1}.$$
\begin{enumerate}
\item
The head of $\m{\alpha,1}$ (i.e. $\m{\alpha,1} / k_1
\m{\alpha,1} $) is isomorphic to the functor  $\kappa(iso_{(x,\alpha)})$ where $iso_{(x,\alpha)}$ is a simple object in $\Gq$.
\item
For $d>0$ 
$$k_d \m{\alpha,1}/k_{d+1} \m{\alpha,1} \simeq L^{d+1}_\alpha $$
where $L^{d+1}_\alpha$ is a simple object of the category $\Fq$ that is neither in the image of
$\iota$ nor in the image of $\kappa$. 

The functor $L^{d+1}_\alpha$ is a subfunctor of $\iota(\Lambda^{d+1}) \otimes \kappa(iso_{(x,\alpha)})$, where $\Lambda^{d+1}$ is the $(d+1)$st exterior power functor.
\end{enumerate}

\end{enumerate}
\end{thm*}

This theorem and the forthcoming paper \cite{Vespa3} lead us to conjecture that there are only three types of simple objects in the category $\Fq$: those in the image of the functor $\iota$, those in the image of the functor $\kappa$ and those which are subfunctors of a tensor product of the form: $\iota(S) \otimes \kappa(T)$ where $S$ is a simple object in $\F$ and $T$ is a simple object in $\Gq$.

This paper is divided into seven sections. Section $1$ recalls the definition of the category $\Fq$ and the results obtained in \cite{math.AT/0606484}. Section $2$ gives a general definition of the mixed functors $\Mn$ as subfunctors of the tensor product $\iota(P^{\F}_{V}) \otimes \kappa(iso_D)$ in
$\Fq$, where $V$ is an object in $\E^f$, $D$ is a quadratic vector space, $\eta$ is an
element in the dual of $V \otimes D$, $P^{\F}_{V}$ is the standard projective object of $\F$ obtained by the Yoneda lemma and $iso_D$ is an isotropic functor in $\Gq$. Section $3$ studies the mixed functors $\Mn$ such that $\mathrm{dim}(D)=1$ and $\mathrm{dim}(V)=1$. We define, in particular, the subfunctor $\m{\alpha,1}$ of the mixed functor $\M{\alpha,1}$ given in the second point of the previous theorem. In section $4$, we deduce a filtration of the functor $\iota(P^{\F}_{\FF}) \otimes  \kappa(iso_{(x,\alpha)})$ from the polynomial filtration in the category $\F$. Section $5$ gives a filtration of the functors $\m{\alpha,1}$, defined in section $3$, and we obtain the existence of a natural map from the subquotients of this filtration to the functors $\iota(\Lambda^{n}) \otimes  \kappa(iso_{(x,\alpha)})$, by relating this filtration to that introduced in the previous section. Section $6$ gives the structure of the functors $\iota(\Lambda^{n}) \otimes  \kappa(iso_{(x,\alpha)})$. We define the functors $L^n_{\alpha}$ and prove their simplicity. Section $7$ proves the structure of the functors $\M{\alpha,1}$ given in the previous theorem.

The results contained in this paper extend results obtained in the author's PhD thesis \cite{Vespa-these}.
The author wishes to thank her PhD supervisor, Lionel Schwartz, for his guidance, as well
as Geoffrey Powell and Aur\'{e}lien Djament for numerous useful discussions and thank Serge Bouc for suggesting that the methods used in the author's thesis should be sufficient to establish the uniseriality of the functors $\m{\alpha,1}$.

%Bien que les foncteurs mixtes d'ordre supérieur soient mal connus, nous
%avons choisi de présenter dans cette thèse la définition de ces
%foncteurs sous leur forme la plus générale, car nous verrons à la
%section \ref{autres-PV} que d'autres foncteurs mixtes apparaissent
%dans les quotients de la filtration des projectifs $P_V$, pour un
%objet $V$ de \Tq\ de dimension supérieure à deux.
\section{The category \Fq}
We recall in this section some definitions and results about the
category $\Fq$ obtained in \cite{math.AT/0606484}.

Let $\Eq$ be the category having as objects finite dimensional $\FF$-vector spaces equipped with a non
degenerate quadratic form and with morphisms linear maps that
preserve the quadratic forms. By the classification of quadratic forms
over the field $\FF$ (see, for instance, \cite{Pfister}) we know that only spaces of even
dimension can be nondegenerate and, for a fixed even dimension, there are  two non-equivalent nondegenerate spaces, which are distinguished by
the Arf invariant. We will denote by $H_0$ (resp. $H_1$) the
nondegenerate quadratic space of dimension two such that $\mathrm{Arf}(H_0)=0$
(resp. $\mathrm{Arf}(H_1)=1$). The orthogonal sum of two nondegenerate quadratic spaces $(V,q_V)$ and $(W,q_W)$ is, by definition, the quadratic space $(V \oplus W,q_{V \oplus W})$ where $q_{V \oplus W}(v,w)=q_V(v)+q_W(w)$. Recall that the spaces $H_0 \bot H_0$
and $H_1 \bot H_1$ are isomorphic. Observe that the morphisms of $\Eq$
are injective linear maps and this category does not admit push-outs
or pullbacks. There exists a pseudo push-out in $\Eq$ that allows us to generalize the construction of the category of
co-spans of B\'{e}nabou  \cite{Benabou} and thus to define the category $\Tq$ in which there
exist retractions.
\begin{defi} \label{1.1}
The category $\Tq$ is the category having as objects those of $\Eq$
and, for $V$ and $W$ objects in $\Tq$, $\mathrm{Hom}_{\Tq}(V,W)$ is
the set of equivalence classes of diagrams in $\Eq$ of the form $ V
\xrightarrow{f} X \xleftarrow{g} W $  for the equivalence relation
generated by the relation $\mathcal{R}$ defined as follows: $ V \xrightarrow{f} X_1 \xleftarrow{g} W \quad  \mathcal{R}\quad   V
\xrightarrow{u} X_2 \xleftarrow{v} W $ if there exists a morphism
$\alpha$ of $\Eq$ such that $\alpha \circ f=u$ and $\alpha \circ g=v$.
The composition is defined using the pseudo push-out. The morphism of
$\mathrm{Hom}_{\Tq}(V,W)$ represented by the diagram $ V
\xrightarrow{f} X \xleftarrow{g} W $ will be denoted by $[ V
\xrightarrow{f} X \xleftarrow{g} W ]$.
\end{defi}

 By definition,
the category $\Fq$ is the category of functors from $\Tq$ to $\E$.
Hence $\Fq$ is
abelian and has enough projective objects. By the Yoneda lemma, for any object $V$ of $\Tq$, the functor $P_V=\FF \lbrack
\mathrm{Hom}_{\Tq}(V,-) \rbrack $ is a projective object and there is
a natural isomorphism: $\mathrm{Hom}_{\Fq}(P_V, F) \simeq F(V)$,
for all objects $F$ of $\Fq$. The set of functors $\{ P_V | V \in
\mathcal{S} \}$, named the standard projective objects in $\Fq$, is a set of
projective generators of $\Fq$, where $\mathcal{S}$ is a set of
representatives of isometry classes of nondegenerate quadratic
spaces.

There is a forgetful functor $\epsilon: \Tq \rightarrow \E^f$ in $\Fq$, defined by $\epsilon(V)=\mathcal{O}(V)$ and 
$$\epsilon([V \xrightarrow{f} W \bot W' \xleftarrow{g} W ])=p_g\circ
\mathcal{O}(f)$$
where $p_g$ is the orthogonal projection from $W \bot W'$ to $W$ and $\mathcal{O}: \Eq \rightarrow \E^f$ is the functor which forgets the quadratic form.  By the fullness of the functor $\epsilon$ and an argument of essential surjectivity, we obtain the following theorem.
\begin{thm} \label{1.2} \cite{math.AT/0606484}
There is a functor $\iota: \F \rightarrow \Fq$, which is exact,
fully faithful and preserves simple objects.
\end{thm}
In order to define another subcategory of $\Fq$, we consider the
category $\Eqd$ having as objects finite dimensional $\FF$-vector spaces
equipped with a (possibly degenerate) quadratic
form and with morphisms injective linear maps that
preserve the quadratic forms. A useful relation between the categories $\Eq$ and $\Eqd$ is given by the following theorem, which can be regarded as Witt's theorem for
degenerate quadratic forms.

\begin{thm} \label{1.3}
Let $V$ be a nondegenerate quadratic space, $D$ and $D'$
subquadratic spaces (possibly degenerate) of $V$ and
$\underline{f}:D \rightarrow D'$ an isometry. Then, there exists an isometry
$f: V \rightarrow V$ such that the following diagram is commutative:
 $$\xymatrix{
V  \ar[r]^f & V  \\
D \ar@{^{(}->}[u] \ar[r]_{\underline{f}}& D'. \ar@{^{(}->}[u] 
}$$
\end{thm}
\begin{proof}
For a proof of this result, we refer the reader to \cite{Bourbaki} \S4, Theorem 1.
\end{proof}
 The category $\Eqd$ admits
pullbacks; consequently the category of spans $\mathrm{Sp}(\Eqd)$ (\cite{Benabou}) is
defined. By definition, the category $\Gq$ is the category of functors
from $\mathrm{Sp}(\Eqd)$ to $\E$. As in the case of the category $\Fq$, the
category $\Gq$ is abelian and has enough projective objects; by the Yoneda lemma, for any object $V$ of $\mathrm{Sp}(\Eqd)$, the functor $Q_V=\FF \lbrack
\mathrm{Hom}_{\mathrm{Sp}(\Eqd)}(V,-) \rbrack $ is a projective object
in $\Gq$. The category $\Gq$ is related to $\Fq$ by the following theorem.
\begin{thm} \label{1.4} \cite{math.AT/0606484}
There is a functor $\kappa: \Gq \rightarrow \Fq$, which is exact,
fully-faithful and preserves simple objects.
\end{thm}

We obtain the classification of the simple objects of the category
$\Gq$ from the following theorem.
\begin{thm} \label{1.5} \cite{math.AT/0606484}
There is a natural equivalence of categories 
$$\Gq \simeq \prod_{V \in \mathcal{S}} \FF[O(V)]-mod$$
where $\mathcal{S}$ is a set of representatives of isometry classes of
quadratic spaces (possibly degenerate) and $O(V)$ is the orthogonal group.
\end{thm}
The object of $\Gq$ that corresponds, by this equivalence, to the
module $\FF[O(V)]$ is the isotropic functor $iso_V$, defined in
\cite{math.AT/0606484}. The family of isotropic functors forms a set of projective generators and injective
cogenerators of $\Gq$. Recall that the isotropic functor $iso_V: \mathrm{Sp}(\Eqd) \rightarrow \E$ of $\Gq$  is the image of $Q_V$ by the morphism $a_V: Q_V \rightarrow DQ_V$ which corresponds by the Yoneda lemma to the element $(\mathrm{Id}_V)^*$ of
$DQ_V(V)$, where $(\mathrm{Id}_V)^*$ is defined by:
$$(\mathrm{Id}_V)^*([\mathrm{Id}_V])=1 \quad \mathrm{and} \quad
(\mathrm{Id}_V)^*([f])=0 \mathrm{\ for\ all\ } f \ne \mathrm{Id}_V$$
where we denote by $[f]$ a canonical generator of $DQ_V(V) \simeq \FF[\mathrm{End}_{\mathrm{Sp}(\Eqd)} (V)]$.
This definition and that of the functor $\kappa: \Gq \rightarrow \Fq$ give rise to the following more concrete definition of the functor $iso_V$ which will be useful below.
\begin{prop} \label{1.6}
The following equivalent definition of the functor $\kappa(iso_V)$ holds.
\begin{itemize}
\item
For $W$ an object of $\Tq$:
$$\kappa(iso_V)(W)=\FF[\mathrm{Hom}_{\Eqd}(V,W)].$$
\item
For a morphism $m=[W \xrightarrow{f} Y \xleftarrow{g} X]$ in $\Tq$ and a canonical generator $[h]$ of $\kappa(iso_V)(W)$, we consider the following diagram in $\Eqd$:
$$\xymatrix{
 &  & X \ar[d]^{g}\\
V \ar[r]_h & W \ar[r]_-f & Y.  }$$
If the pullback of this diagram in $\Eqd$ is $V$, this gives rise to a unique morphism $h': V \rightarrow X$  in $\Eqd$, such that $f \circ h=g \circ h'.$
In this case, $\kappa(iso_V)(m)[h]=[h'].$
Otherwise, $\kappa(iso_V)(m)[h]=0.$
\end{itemize}
\end{prop}

\begin{nota}
In this paper, a canonical generator of $\kappa(iso_D)(W)$ will be denoted by: $[D \xrightarrow{h} W]$ or, more simply, by $[h]$.
\end{nota}

We end this section by a useful corollary of Theorem \ref{1.4} and
\ref{1.5}. For $\alpha \in \{ 0,1 \}$, let $(x,\alpha)$ be the degenerate quadratic space of
dimension one generated by $x$ such that $q(x)=\alpha$. 
\begin{cor} \label{1.8} 
The functors $\kappa(iso_{(x,0)})$ and $\kappa(iso_{(x,1)})$
are simple in $\Fq$.
\end{cor}
\begin{proof}
It is a straightforward consequence of the triviality of the orthogonal groups $O(x,0)$ and $O(x,1)$.
\end{proof}

%%%%%%%%%%%%%%%%%%%%%%%%%%%%%%%%%%%%%%%%%%%%%%%

\section{Definition of the mixed functors}
The aim of this section is to define the mixed functors: for this, we
consider the functors $\iota(P^{\F}_{V}) \otimes \kappa(iso_D)$ in
$\Fq$ where $V$ is an object in $\E^f$, $P^{\F}_{V}$ is the standard projective object of $\F$ obtained by the Yoneda lemma, $D$ is an object in $\Eqd$,  and $\iota: \F \rightarrow \Fq$
and $\kappa: \Gq \rightarrow \Fq$ are the functors defined in Theorem \ref{1.2} and \ref{1.4} respectively. A canonical generator of $P^{\F}_{V}(W) \simeq \FF [ \mathrm{Hom}_{\E^f}(V,W)]$ will be denoted by $[f]$.

\begin{nota}
In this paper, the bilinear form associated to a quadratic space $V$
will be denoted by $B_V$.
\end{nota}

\begin{prop} \label{2.2}
Let $D$ be an object in $\Eqd$, $V$ be an object in $\E^f$, $\eta$ be an
element in the dual of $V \otimes D$ and $W$ be an object in $\Tq$, the subvector space of
$(\iota(P^{\F}_{V}) \otimes \kappa(iso_D))(W)$ generated by the
elements:
$$[f] \otimes [ D \xrightarrow{h} W]$$
such that 
$$\forall v \in V, \forall d \in D \quad B_W(f(v), h(d))=\eta(v \otimes d)$$
defines a subfunctor of $\iota(P^{\F}_{V}) \otimes \kappa(iso_D)$
which will be denoted by $\Mn$ and called the mixed functor associated
to $V$, $D$ and $\eta$.
\end{prop}

\begin{proof}
It is sufficient to verify that, for each morphism $M=[W \xrightarrow{k} Y \xleftarrow{l} Z]$ of $\Tq$ and each generator $[f] \otimes [ D \xrightarrow{h} W]$ of $\Mn(W)$, 
$$\Mn(M)([f] \otimes [D \xrightarrow{h} W]) \in  \Mn(Z).$$
Consider the following diagram in $\Eqd$:
$$\xymatrix{
 &  & Z \ar[d]^{l}\\
D \ar[r]_h & W \ar[r]_-k & Y.  }$$
\begin{itemize}
\item
If the pullback of this diagram in $\Eqd$ is $D$, namely if $k \circ
h(D) \subset l(Z)$, this gives rise to a unique morphism $h'$, from $D$ to
$Z$ in $\Eqd$, such that $k \circ h=l \circ h'$ that is, the following diagram commutes
$$\xymatrix{
D \ar[rr]^{h'} \ar[d]_{\mathrm{Id}} & &\ar[d]^{l}  Z\\
D \ar[r]_h & W \ar[r]_-k & Y. }$$
In this case, by Proposition \ref{1.6}, we have: 
$$\Mn(M)([f] \otimes [ D \xrightarrow{h} W])= ([p_l \circ k \circ f] \otimes [ D \xrightarrow{h'} Z]) $$
where $p_l$ is the orthogonal projection associated to $l$. For an element $v$ in $V$ and $d$ in $D$, we have:
$$B_V(f(v),h(d))=B_Y(k \circ f(v), k \circ h(d)).$$
Since the pullback of the diagram considered previously  is $D$, we have: $k \circ h(D) \subset l(Z).$
Consequently:
$$B_V(f(v),h(d))=B_Z(p_l \circ k \circ f(v), p_l \circ k \circ h(d))=B_Z(p_l \circ k \circ f(v), h'(d)).$$
Thus, if $B_V(f(-),h(-))=\eta$ then $B_Z(p_l \circ k \circ f(-), h'(-))=\eta$, therefore
$([p_l \circ k \circ f] \otimes [D \xrightarrow{h'} Z])$ is an element of $\Mn(Z)$.
\item
Otherwise we have, by Proposition \ref{1.6}: 
$$ \Mn(M)([f] \otimes [D \xrightarrow{h} W])= 0.$$
\end{itemize}

\end{proof}

\begin{rem}
The terminology "mixed functors" is chosen to reflect the fact that these
functors are subfunctors of a tensor product of a functor coming
from the category $\F$ and a functor coming from the category $\Gq$.
\end{rem}
We obtain the following decomposition of the functors
$\iota(P^{\F}_{V}) \otimes \kappa(iso_D)$.
\begin{lm} \label{2.4}
For $D$ an object in $\Eqd$ and $V$ an object in $\E^f$ we have:
$$\iota(P^{\F}_{V}) \otimes \kappa(iso_D)= \underset{\eta \in (V \otimes D)^*}{\oplus} \Mn.$$
\end{lm} 
\begin{proof}
For two different elements $\eta$ and $\eta'$ in $(V \otimes D)^*$, we have:
$$\mathrm{Mix}_{V, D,\eta}(W) \cap \mathrm{Mix}_{V, D,\eta'}(W)=\{ 0\}$$
for $W$ an object in $\Tq$. Thus, we have the decompositions: $(\iota(P^{\F}_{V}) \otimes \kappa(iso_D))(W)= (\underset{\eta \in
  (V \otimes D)^*}{\oplus} \Mn)(W)$, for all objects \mbox{$W$ in~$\Tq$}. Since $\Mn$ is a subfunctor of $\iota(P^{\F}_{V}) \otimes \kappa(iso_D)$ by Proposition \ref{2.2}, we deduce the result.
\end{proof}
%\begin{rem}
%We will give a more general definition of mixed functors in a future
%paper, by using the description of the morphisms of $\Tq$ in terms
%of morphisms of $\mathrm{Sp}(\Eqd)$ and $\E$. The functors define previously will appear as a subfamily of
%these generalised mixed functors.
%\end{rem}

\begin{rem} \label{2.5}
In the definition of the mixed functors, we don't impose that $h(D) \cap f(V)=\{ 0\}$. Nevertheless, we can define similar functors with this condition, which give rise to quotient functors to the mixed functors defined in Proposition \ref{2.2}. These functors will be useful for a later general study of the mixed functors.
\end{rem}

%%%%%%%%%%%%%%%%%%%%%%%%%%%%%%%%%%%%%%%%%%%%%%%%%%%%%

\section{The functors \Mn\ such that $\mathrm{dim}(D)=1$ and $\mathrm{dim}(V)=1$}
The aim of this section is to give some general results about the four simplest mixed functors of $\Fq$ obtained in the case of $\mathrm{dim}(D)=\mathrm{dim}(V)=1$. The motivation of the particular interest in this case is the study of the projective generators $P_{H_0}$ and $P_{H_1}$ of $\Fq$. In fact, we prove in \cite{Vespa3}, that the mixed functors that are direct summands of these two standard  projective generators of $\Fq$ verify the conditions $\mathrm{dim}(D)=\mathrm{dim}(V)=1$.

When $V$ and $D$ are spaces of dimension one, we will denote by  $\Mab$, where $\alpha$ and $\beta$  are elements of $\{0,1\}$, the functor \Mn\ such that $V \simeq \FF$, $D \simeq (x,\alpha)$ and $\eta=\beta$.
We have the following result.
\begin{lm} \label{3.1}
Let $W$ be an object in $\Tq$, if $[f] \otimes [(x,\alpha) \xrightarrow{h} W]$ is a canonical generator of $\Mab(W)$, then $[f+h(x)] \otimes [(x,\alpha) \xrightarrow{h} W]$ is a canonical generator of $\Mab(W)$.
\end{lm}
\begin{proof}
This is a straightforward consequence of the fact that the bilinear form associated to a quadratic form is alternating.
\end{proof}
In order to make this symmetry clearer in the set of canonical generators of $\Mab(W)$ and to introduce an action of the symmetric group $\mathfrak{S}_2$ on this set, we use a slightly different description of the canonical generators of $\Mab(W)$ corresponding to a reindexing of these canonical generators.

\begin{defi} \label{3.2}
For $\alpha$ and $\beta$ elements of $\{ 0,1 \}$, we consider the following set:
$$N^W_{\alpha,\beta}=\{(w_1,w_2)\  |\  w_1 \in W, w_2 \in W,\ q(w_1+w_2)=
\alpha,\ B(w_1,w_2)=\beta\}.$$
\end{defi}

We have the following result.
\begin{lm} \label{3.3}
For $D \simeq (x,\alpha)$ and $\eta=\beta$, we have:
$$\Mab(W) \simeq \FF[ N^W_{\alpha,\beta} ]$$
where $W$ is an object in $\Tq$.
\end{lm}
\begin{proof}
The generator of the vector space $V$ of dimension one will be denoted by $a$. There is an isomorphism:
$$\begin{array}{llll}
f_W:& \Mab(W) & \rightarrow & \FF[N^W_{\alpha,\beta} ]\\
    &     [f] \otimes [(x,\alpha) \xrightarrow{h} W]  & \mapsto    &  [(f(a)+h(x), f(a)],
    \end{array}$$
of which the inverse is given by:
$$\begin{array}{llll}
f_W^{-1}:& \FF[N^W_{\alpha,\beta} ] & \rightarrow & \Mab(W) \\
    & [(w_1,w_2)]      & \mapsto    & [k] \otimes [(x, \alpha) \xrightarrow{l} W]
\end{array}$$
where $k:V \rightarrow W$ is defined by $k(a)=w_2$ and $l:(x, \alpha) \rightarrow W$ is defined by $l(x)=w_1+w_2$.
\end{proof}

\begin{nota} \label{3.4}
Henceforth, we will use the basis given by the set $N^W_{\alpha,\beta} $ to represent the elements of $\Mab(W)$. 
\end{nota}
Thus, the canonical generator $[f] \otimes [(x,\alpha) \xrightarrow{h} W] $ of $\Mab(W)$ is represented by $[(f(a)+h(x), f(a))]$  and $[f +h(x)] \otimes [(x,\alpha) \xrightarrow{h} W] $, which is also a canonical generator of $\Mab(W)$ by Lemma \ref{3.1}, is represented by $[(f(a), f(a)+h(x))]$.

We have the following lemma.
\begin{lm} \label{3.5}
The symmetric group $\mathfrak{S}_2$ acts on the functor $\M{\alpha,\beta}$.
\end{lm}
\begin{proof}
Let $W$ be an object of $\Tq$, we define an action of $\mathfrak{S}_2=\{\mathrm{Id},
\tau \}$ on $\M{\alpha,\beta}(W)$ by:
$$\tau.[(w_1,w_2)]=[(w_{2},w_1)].$$
We leave the reader to verify that the linear maps $\tau_W$
$$\begin{array}{llll}
\tau_W:& \Mab(W) &\rightarrow& \Mab(W)\\
         & [(w_1,w_2)] & \mapsto &  [(w_2,w_1)]
\end{array}$$
define a natural transformation.
\end{proof}

This lemma allows us to define an object in $\Fq$ by considering the invariants by this action.

\begin{defi} \label{3.6}
Let $\mab$ be the subfunctor of $\Mab$ defined by considering the invariants of $\Mab(W)$ by the action of the symmetric group $\mathfrak{S}_2$.
%$$\begin{array}{llll}
%\mab: & \Tq & \rightarrow & \E \\
 %     &  W & \longmapsto & (\Mab(W))^{\sigma_2}.
%\end{array}$$
\end{defi}

In the following, we will focus on study  the functors \M{0,1}\ and \M{1,1}. These two functors are particularly interesting since they are direct summands of $P_{H_0}$ and $P_{H_1}$ (see \cite{Vespa3}). 

We have the following lemma.
\begin{lm} \label{3.7}
Let $W$ be an object in $\Tq$ and $[(w_1,w_2)]$ be a generator of $\M{\alpha,1}(W)$, then the vectors $w_1$ and $w_2$ are linearly independent.
\end{lm}

\begin{proof}
This is a straightforward consequence of the fact that the bilinear form $B$ is alternating.
\end{proof}

We deduce the following lemma.
\begin{lm} \label{3.8}
Let $W$ be an object in $\Tq$, the action of $\mathfrak{S}_2$ on the set of canonical generators of $\M{\alpha,1}(W)$ is free.
\end{lm}

\begin{proof}
For a canonical generator $[(w_1,w_2)]$ of $\M{\alpha,1}(W)$, since the vectors $w_1$ and $w_2$ are linearly independent by Lemma \ref{3.7}, we have $w_1 \ne w_2$. Hence, the action of $\mathfrak{S}_2$ is free.
\end{proof}

\begin{rem} \label{3.9}
We deduce from Lemma \ref{3.7} that the two functors $\M{\alpha,1}$, coincide with the functors mentioned in Remark \ref{2.5}.
\end{rem}

We give the following general result about the free actions of the group $\mathfrak{S}_2$.
\begin{lm} \label{3.10}
If $A$ is a finite set equipped with a free action of the group $\mathfrak{S}_2$ then,
there exists a short exact sequence of $\mathfrak{S}_2$-modules:
$$0 \rightarrow \FF[A]^{\mathfrak{S}_2} \rightarrow \FF[A] \rightarrow \FF[A]^{\mathfrak{S}_2} \rightarrow
0.$$
\end{lm}
\begin{proof}
We deduce from the action of $\mathfrak{S}_2$ on $A$, the existence of the canonical inclusion of the invariants in $\FF[A]$: $\xymatrix{
\FF[A]^{\mathfrak{S}_2} \ar@{^{(}->}[r]^{f}& \FF[A].
}$
The norm $\FF[A] \xrightarrow{1+\tau} \FF[A]$ induces a linear map: $\FF[A]  \xrightarrow{g} F[A]^{\mathfrak{S}_2}$, 
such that the composition 
$$\xymatrix{
\FF[A]^{\mathfrak{S}_2} \ar@{^{(}->}[r]^{f}& \FF[A]  \ar[r]^g & \FF[A]^{\mathfrak{S}_2}
}$$
is trivial. We verify that this defines a short exact sequence.
\end{proof}
We deduce the following proposition.
\begin{prop} \label{3.11}
There exists a short exact sequence
\begin{equation} \label{Mix2-sec}
0 \rightarrow \m{\alpha,1} \rightarrow \M{\alpha,1} \rightarrow \m{\alpha,1} \rightarrow 0.
\end{equation}
\end{prop}
\begin{proof}
This is a straightforward consequence of Lemma \ref{3.8} and Lemma \ref{3.10}.
\end{proof} \label{3.12}
\begin{nota}
We will denote by $[\{ w_1, w_2 \}]$ the image of the element  $[(w_1, w_2)]$ of $ \M{\alpha,1}(W)$ in $ \m{\alpha,1}(W)$ by the surjection 
$\xymatrix{
\M{\alpha,1}(W) \ar@{->>}[r] & \m{\alpha,1}(W).
}$
\end{nota}

\begin{rem}
Lemma \ref{3.8} has no analogue for the functors \M{0,0}\ and \M{1,0} since, in these two cases, the action of the group $\mathfrak{S}_2$ is not free. Nevertheless, we can apply similar arguments to the functors mentioned in Remark \ref{2.5} corresponding to the functors \M{0,0}\ and \M{1,0}. 
\end{rem}

\begin{rem}
It is shown in \cite{Vespa3} that the functors $\M{0,1}$ and $\M{1,1}$ are indecomposable. Consequently, the short exact sequence (\ref{Mix2-sec}) is not split for the functors $\M{0,1}$ and $\M{1,1}$.
\end{rem}

%-----------------------------------------------------------------------
\section{Study of the functor $\iota(P^{\F}_{\FF}) \otimes  \kappa(iso_{(x,\alpha)})$} 

By Proposition \ref{2.2}, the functor $\M{\alpha,1}$ is a subfunctor of $\iota(P^{\F}_{\FF}) \otimes  \kappa(iso_{(x,\alpha)})$. Consequently, in order to obtain the decomposition of $\M{\alpha,1}$, we study, in this section, the functor $\iota(P^{\F}_{\FF}) \otimes  \kappa(iso_{(x,\alpha)})$.

\subsection{Filtration of the functors $\iota(P^{\F}_{\FF}) \otimes  \kappa(iso_{(x,\alpha)})$}
We define, below, the filtration of the functors $\iota(P^{\F}_{\FF}) \otimes  \kappa(iso_{(x,\alpha)})$ induced by the polynomial filtration of the functor $P^{\F}_{\FF}$ in the category $\F$. First we recall  the essential results concerning the polynomial functors in the category $\F$. We refer the interested reader to \cite{HLS1}, \cite{Sch} and \cite{KuhnI} for details on the subject.

\begin{nota}
Henceforth, in order to simplify the notation, we will denote the functor $\iota(P^{\F}_{\FF}) \otimes  \kappa(iso_{(x,\alpha)})$ by $P_{\mathbb{F}} \otimes  iso_{\alpha}$ and, if $F \neq P^{\F}_{\FF}$, we will denote the functor $\iota(F) \otimes  \kappa(iso_{(x,\alpha)})$ by $F \otimes  iso_{\alpha}$.
\end{nota}

\begin{defi}
Let $F$ be an object in $\F$ and $d$ an integer, the functor $q_d F$ is the largest polynomial quotient of degree $d$ of the functor $F$.
\end{defi}
\begin{nota}
We denote by $k_d F$ the kernel of 
$\xymatrix{
F \ar@{->>}[r]& q_d F}$.
\end{nota}
We have the following result.
\begin{prop} \label{4.4}
For an integer $d$, the functors $k_d F$ define a decreasing filtration of the functor $F$.
\end{prop}

Thus, for the standard projective functor $P_{\mathbb{F}}$, we have the following short exact sequence:
\begin{equation} \label{filt-sec1}
0 \rightarrow k_d P_{\mathbb{F}} \rightarrow P_{\mathbb{F}}
\xrightarrow{f_d} q_d P_{\mathbb{F}} \rightarrow 0.
\end{equation}
Furthermore, the decreasing filtration of $P_{\mathbb{F}}$ given by the functors $k_d P_{\mathbb{F}} $ is separated (that is $\cap k_d P_{\mathbb{F}} =0$).

We recall below the description of the vector space $k_d
P_{\mathbb{F}}(V)$ for $V$ an \mbox{object in $\mathcal{E}^f$.}
\begin{prop} \label{4.5}
The vector space $k_d P_{\mathbb{F}}(V)$ is generated by the elements 
$$\sum_{z \in \mathcal{L}} [z]$$
where $\mathcal{L}$ is a subvector space of $V$ of dimension $d+1$.
\end{prop}
\begin{nota}
The subvector space of $V$ or subquadratic space of $(V, q_V)$ generated by $v_1, \ldots v_n$ will be denoted by $\mathrm{Vect}( v_1, \ldots v_n )$. 
\end{nota}
The subquotients of the filtration of the functor $P_{\mathbb{F}}$ are given in the following proposition.
\begin{prop} \label{4.6}
For $d$ a non-negative integer, there exists a short exact sequence
\begin{equation} \label{filt-sec2}
0 \rightarrow k_{d+1} P_{\mathbb{F}} \rightarrow k_d P_{\mathbb{F}}
\xrightarrow{g_d} \Lambda^{d+1} \rightarrow 0
\end{equation}
where $ \Lambda^{d+1} $ is the $(d+1)$-th exterior power and the map $g_d$ is defined in the following way: for $V$ an
object in $\E^f$ and $\mathcal{L}$ the subvector space of $V$ of dimension $d+1$ generated by the elements $l_1, \ldots, l_{d+1}$ of $V$, we have:
$$(g_d)_V(\sum_{z \in \mathcal{L}} [z])=l_1 \wedge \ldots \wedge l_{d+1}.$$

\end{prop}
By taking tensor product with the isotropic functors $iso_{\alpha}$, Proposition \ref{4.4} and the short exact sequences (\ref{filt-sec1}) and (\ref{filt-sec2}) give rise to the following result.
\begin{cor} \label{4.7}
\begin{enumerate}
\item
For $d$ an integer, the functors $(k_d P_{\mathbb{F}})
\otimes iso_{\alpha}$ define a decreasing separated filtration of the functor  $ P_{\mathbb{F}} \otimes iso_{\alpha}$.
\item
There exist the following short exact sequences in \Fq:
\begin{equation} \label{filt-sec3}
0 \rightarrow (k_d P_{\mathbb{F}}) \otimes iso_{\alpha} \rightarrow P_{\mathbb{F}} \otimes iso_{\alpha}
\xrightarrow{f_d \otimes iso_{\alpha}} (q_d P_{\mathbb{F}}) \otimes iso_{\alpha} \rightarrow 0;
\end{equation}

\begin{equation} \label{filt-sec4}
0 \rightarrow (k_{d+1} P_{\mathbb{F}}) \otimes iso_{\alpha} \rightarrow (k_d
P_{\mathbb{F}}) \otimes iso_{\alpha}
\xrightarrow{g_d \otimes iso_{\alpha}} (\Lambda^{d+1}) \otimes iso_{\alpha} \rightarrow 0.
\end{equation}
\end{enumerate}
\end{cor}
\begin{rem}
We obtain similar results by taking the tensor product between the short exact sequence \ref{filt-sec2} and the isotropic functors $iso_D$. This will be useful for a general study of the mixed functors.
\end{rem}
%%%%%%%%%%%%%%%%%%%%%%%%%%%%%%%%%%%%%%%%%

\section{Filtration of the functors  $\m{\alpha,1}$} \label{}

In this section, we define a filtration of the functors $\m{\alpha,1}$ that we will relate below to the filtration of $P_{\mathbb{F}}
\otimes iso_{\alpha}$ obtained in the previous section.

\begin{defi}
For $V$ an object of $\Tq$ and $d$ an integer, let $k_d
\m{\alpha,1}(V)$ be the subvector space of $\m{\alpha,1}(V)$ generated by the elements: 
$$\sum_{z \in \mathcal{L}} [ \{ x+z,y+z \} ],$$
where $[\{ x,y \} ] \in \m{\alpha,1}(V)$ and $\mathcal{L}$ is a subvector space of $\mathrm{Vect}( x+y )^{\bot}$ of dimension $d$.
\end{defi}

\begin{prop} \label{5.2}
The spaces $k_d \m{\alpha,1}(V)$, for $V$ an object in $\Tq$, define a subfunctor of $\m{\alpha,1}$.
\end{prop}

\begin{proof}
For a morphism $T$ in $\mathrm{Hom}_{\Tq}(V,W)$, it is straightforward to check, by definition of $\M{\alpha,1}(T)$, that the image of $k_d \m{\alpha,1}(V)$ by
$$\m{\alpha,1}(T): k_d \m{\alpha,1}(V) \rightarrow  \m{\alpha,1}(W)$$
is a subvector space of $k_d \m{\alpha,1}(W).$
\end{proof}

\begin{prop}
The functors $k_d \m{\alpha,1}$ define a separated decreasing filtration of the functor $\m{\alpha,1}$
$$\ldots \subset k_d \m{\alpha,1} \subset \ldots \subset k_1 \m{\alpha,1} \subset k_0 \m{\alpha,1}= \m{\alpha,1}.$$
\end{prop}

\begin{proof}
To show that the filtration is decreasing, it is sufficient to prove, for an object $V$ of $\Tq$, that there is an inclusion of vector spaces:
$$k_{d+1} \m{\alpha,1}(V) \subset k_d \m{\alpha,1}(V).$$
We consider a generator of $k_{d+1} \m{\alpha,1}(V)$: $ v=\sum_{z \in \mathcal{L}} [ \{ x+z,y+z \} ]$
where  $[\{ x,y \} ] \in \m{\alpha,1}(V)$ and $\mathcal{L}$ is the subvector space of $\mathrm{Vect}( x+y )^{\bot}$ of dimension $d+1$ generated by the elements $l_1,
\ldots, l_{d+1}$, and we consider the following decomposition into direct summands: $\mathcal{L}= \mathrm{Vect}( l_1, \ldots, l_d ) \oplus \mathrm{Vect}(l_{d+1})=\mathcal{L}'\oplus \mathrm{Vect}( l_{d+1}).$
By considering separately the elements $z$ in $\mathcal{L}$ with an nonzero component on $l_{d+1}$ and those with a zero component on $l_{d+1}$, we obtain:
$$v=\sum_{z \in \mathcal{L}'} [ \{ x+z,y+z \} ] +\sum_{z \in
  \mathcal{L}'} [ \{ x+l_{d+1}+z,y+l_{d+1}+z \} ].$$
 By definition of $k_d  \m{\alpha,1}$ we have:  $\sum_{z \in \mathcal{L}'} [ \{ x+z,y+z \} ] \in  k_d
 \m{\alpha,1}(V)$
  and   $ \sum_{z \in
  \mathcal{L}'} [ \{ x+l_{d+1}+z,y+l_{d+1}+z \} ] \in  k_d
  \m{\alpha,1}(V)$
  since $[ \{ x+l_{d+1},y+l_{d+1} \} ] \in  \m{\alpha,1}(V).$ Consequently $v  \in  \m{\alpha,1}(V).$

  One verifies easily that the filtration is separated.                 
\end{proof}

In the following result, we relate the filtration of the functors $\m{\alpha,1}$ and the filtration of $P_{\mathbb{F}}
\otimes iso_{\alpha}$ obtained from the polynomial filtration in Corollary \ref{4.7}.
\begin{lm} \label{5.4}
The composition
$$\xymatrix{
k_d \m{\alpha,1} \ar@{^{(}->}[r] & \m{\alpha,1}  \ar@{^{(}->}[r] & \M{\alpha,1}
  \ar@{^{(}->}[r] & P_{\mathbb{F}}
\otimes iso_{\alpha} \ar[r]^{f_d \otimes iso_{\alpha} } & q_d P_{\mathbb{F}}
\otimes iso_{\alpha}
}$$
 is zero.
\end{lm}
\begin{proof}
Let $V$ be an object in $\Tq$ and $v$ a generator of $k_d
\m{\alpha,1}(V)$:  $v=\sum_{z \in \mathcal{L}} [\{ y+z, y'+z \} ]$
where $[\{ y,y' \} ] \in  \m{\alpha,1}(V)$ and $\mathcal{L}$ is the subvector space of $\mathrm{Vect}( y+y' )^{\bot}$ of dimension $d$ generated by the elements $l_1, \ldots l_d$.

Let $h$ be the element of $\mathrm{Hom}_{\Eq^d}((x,\alpha),V)$ such that $h(x)=y+y'$, we have: 
$$\begin{array}{l}
\sum_{z \in \mathcal{L}} ( [y+z] + [y'+z] ) \otimes
[h]\\
            =\sum_{z \in \mathcal{L}} ( [y+z]+[z] + [y'+z]+[z] ) \otimes
[h]\\
            =\sum_{z \in \mathcal{L}} ( [y+z]+[z]) \otimes [h] + \sum_{z \in \mathcal{L}} ( [y'+z]+[z]) \otimes [h]\\
           =\sum_{z' \in \mathcal{L} \oplus \mathrm{Vect}( y)} [z']
           \otimes [h] + \sum_{z'' \in \mathcal{L} \oplus \mathrm{Vect}( y') } [z''] \otimes [h].
\end{array}$$
By Proposition \ref{4.5} we have:
$$\sum_{z' \in \mathcal{L} \oplus \mathrm{Vect}( y )} [z'] \in  k_d
P_{\mathbb{F}}(V) \mathrm{\quad and\quad } \sum_{z'' \in \mathcal{L} \oplus \mathrm{Vect}( y') } [z''] \in  k_d
P_{\mathbb{F}}(V),$$
hence $f_d \otimes iso_{\alpha} (v)=0.$
\end{proof}
We have the following result.
\begin{prop} \label{5.5}
\begin{enumerate}
\item
There exists a monomorphism $i_d$ from $k_d \m{\alpha,1}$ to $k_d P_{\mathbb{F}} \otimes iso_{\alpha}$.
\item
There exists a natural map:
$$k_d \m{\alpha,1}/k_{d+1} \m{\alpha,1} \rightarrow \Lambda^{d+1} \otimes
iso_{\alpha}$$
induced by $i_d: k_d \m{\alpha,1} \rightarrow k_d P_{\mathbb{F}} \otimes iso_{\alpha}$.
\end{enumerate}
\end{prop}
\begin{proof}
The first point is a direct consequence of Lemma \ref{5.4} and the short exact sequence (\ref{filt-sec3}).

We deduce the second point from the following commutative diagram given by the first point and from the short exact sequence (\ref{filt-sec4}).
\begin{equation} \label{filt-diagramme}
\xymatrix{
0 \ar[r] & k_{d+1} \m{\alpha,1} \ar[r] \ar@{^{(}->}[d]_{i_{d+1}} &  k_{d}
\m{\alpha,1} \ar[r]  \ar@{^{(}->}[d]_{i_d} & (k_d \m{\alpha,1} / k_{d+1}
\m{\alpha,1}) \ar[r]  \ar@{->}[d] & 0\\
0 \ar[r] & k_{d+1} P_{\mathbb{F}} \otimes iso_{\alpha} \ar[r] &
k_{d} P_{\mathbb{F}} \otimes iso_{\alpha} \ar[r]_{g_d \otimes iso_{\alpha}}  &  \Lambda^{d+1} \otimes
iso_{\alpha} \ar[r] & 0
}
\end{equation}
\end{proof}
Hence, to obtain the composition factors of the functors $\M{\alpha,1}$ we study the functors $\Lambda^n \otimes iso_{\alpha}$ in the following section.

\begin{rem}
As a consequence of Theorem \ref{7.1} we obtain that the natural map $k_d \m{\alpha,1}/k_{d+1} \m{\alpha,1} \rightarrow \Lambda^{d+1} \otimes
iso_{\alpha}$ is a monomorphism.
\end{rem}

To conclude this section, we prove that the functors $\M{\alpha,1}$ are infinite. For this, we need the following lemma.
\begin{lm} \label{5.7}
Let $V$ be an object in $\Tq$ of dimension greater than $d+1$, such that
$\m{\alpha,1}(V) \ne \{0\}$, $[\{ y,y'
\}]$ be a canonical generator of $\m{\alpha,1}(V)$ and $v_1, \ldots
v_d$, be $d$ linearly independent elements in $\mathrm{Vect}( y+y')^{\bot}$, then:
$$(g_d \otimes iso_{\alpha}) \circ i_d(\sum_{z \in \mathrm{Vect}( v_1,
  \ldots, v_{d} ) } [\{ y+z,y'+z \}]) =  v_1 \wedge \ldots \wedge v_d \wedge
                                (y+y') \otimes [h]  $$
where $i_d$ is the monomorphism from  $k_d
\m{\alpha,1}$ to $k_d P_{\mathbb{F}} \otimes iso_{\alpha}$ defined in the first point of Proposition \ref{5.5} and $h$ the element of $\mathrm{Hom}_{\Eq}((x,\alpha), V)$ such that $h(x)=y+y'$.
\end{lm}

\begin{proof}
We have, by definition of $i_d$ and $g_d \otimes iso_{\alpha}$,
$$\begin{array}{ll}
(g_d \otimes iso_{\alpha}) \circ i_d(\sum_{z \in \mathrm{Vect}( v_1,
  \ldots, v_{d} )} [\{ y+z,y'+z \}]) \\
  = (g_d \otimes
iso_{\alpha})(\sum_{z  \in \mathrm{Vect}( v_1, \ldots, v_{d} ) }  ([y+z]+[y'+z])\otimes [h])\\
                                =(v_1 \wedge \ldots \wedge v_d \wedge
                                y + v_1 \wedge \ldots \wedge v_d \wedge
                                y' ) \otimes [h]\\
                                = v_1 \wedge \ldots \wedge v_d \wedge
                                (y+y') \otimes [h] 

\end{array}$$
\end{proof}
We deduce the following result.
\begin{prop} \label{5.8}
The functors $\M{\alpha,1}$ are infinite.
\end{prop}
\begin{proof}
It is sufficient to prove that the quotients of the filtration of $\m{\alpha,1}$ are nonzero.
For an object $V$ in $\Tq$ of dimension greater than $d$, the space $\m{1,1}(H_0 \bot V)$ contains the nonzero element  $[\{a_0,b_0 \}]$.
Hence, the element of $k_{d} \m{1,1}(H_0 \bot V)$
$$x=\sum_{z \in \mathrm{Vect}( v_1, \ldots, v_{d} ) } [\{ a_0+z,b_0+z
\}] $$
verifies
$$(g_d \otimes iso_{\alpha}) \circ i_d(x)= v_1 \wedge \ldots \wedge v_d \wedge
                                (a_0+b_0) \otimes [h] \ne 0$$
                                                                by Lemma \ref{5.7}. 
Consequently, $(g_d \otimes iso_{\alpha}) \circ i_d(k_d \m{1,1})
                                \ne \{0 \}$ and, by the commutativity of the diagram (\ref{filt-diagramme}) given in the proof of the second point of Proposition \ref{5.5}, we have \mbox{$k_d \m{1,1} /
k_{d+1} \m{1,1}\ne \{ 0\}$.}

In the same way, by considering the element  
$$\sum_{z \in \mathrm{Vect}( v_1, \ldots, v_{d} )} [\{ a_0+z,a_0+b_0+z
\}] \in k_{d} \m{0,1}(H_0 \bot V),$$
we show: $k_d \m{0,1} /
k_{d+1} \m{0,1} \ne \{ 0\}$.
\end{proof}

%-------------------------------------------------------------------------
\section{Structure of the functors $\Lambda^n \otimes iso_{\alpha}$}

By the second point of Proposition \ref{5.5}, there exists a natural map from the subquotients $k_d \m{\alpha,1}/k_{d+1} \m{\alpha,1}$ of the filtration of the functor $\m{\alpha,1}$ by the functors $k_d \m{\alpha,1}$ to the functors $\Lambda^{d+1} \otimes
iso_{\alpha}$. The aim of this section is to study the functors $\Lambda^n \otimes \mathrm{Iso}_{\alpha}$ in order to obtain the composition factors of the functors $\m{\alpha,1}$. It is divided into two subsections: the first concerns the decompositions of the functors
  $\Lambda^n \otimes \mathrm{Iso}_{\alpha}$ by the functors denoted by $L_{\alpha}^n$, and the second concerns the simplicity of the functors
  $L_{\alpha}^n$.

\subsection{Decomposition}
In order  to identify the composition factors of the functor
$\Lambda^n \otimes \mathrm{Iso}_{\alpha}$, we define the following morphisms of
\Fq.

\begin{lm} \label{6.1}
\begin{enumerate}
\item
For an object $V$ in $\Tq$, the linear maps
$$\mu_V:iso_{\alpha}(V) \rightarrow (\Lambda^1 \otimes iso_{\alpha})(V)$$
defined by:
$$\mu_V([(x,\alpha) \xrightarrow{h} V])=h(x) \otimes [(x,\alpha) \xrightarrow{h} V]$$
for a canonical generator $[(x,\alpha) \xrightarrow{h} V]$ in $iso_{\alpha}(V)$, give rise to a monomorphism
$\mu: iso_{\alpha} \rightarrow \Lambda^1 \otimes iso_{\alpha}$ of \Fq.

\item
For an object $V$ in $\Tq$, the linear maps
$$\nu_V:(\Lambda^1 \otimes iso_{\alpha})(V) \rightarrow iso_{\alpha}(V)$$
defined by:
$$\nu_V(w \otimes [(x,\alpha) \xrightarrow{h} V])=B(w,h(x)) [(x,\alpha) \xrightarrow{h} V]$$
for a canonical generator $[(x,\alpha) \xrightarrow{h} V]$ in $iso_{\alpha}(V)$, give rise to an epimorphism
$\nu: \Lambda^1 \otimes iso_{\alpha} \rightarrow iso_{\alpha} $ of \Fq.
\end{enumerate}
\end{lm}
\begin{proof}
\begin{enumerate}
\item
We check the commutativity of the following diagram, for a morphism $T=[V \xrightarrow{f} X \xleftarrow{g} W ]$ of $\mathrm{Hom}_{\Tq}(V,W)$:
$$\xymatrix{
iso_{\alpha}(V) \ar[r]^-{\mu_V} \ar[d]_{iso_{\alpha}(T)}& (\Lambda^1 \otimes iso_{\alpha})(V) \ar[d]^{(\Lambda \otimes iso_{\alpha})(T)}       \\
iso_{\alpha}(W) \ar[r]_-{\mu_W} &  (\Lambda^1 \otimes iso_{\alpha})(W).
}$$
For a canonical generator $[(x,\alpha) \xrightarrow{h} V]$ of $iso_{\alpha}(V)$, we denote by $P$ the pullback of: $(x, \alpha) \xrightarrow{f \circ h} X \xleftarrow{g} W$.
If $P = (x,\alpha)$, we denote by $h'$ the morphism making the following diagram commutative:
$$\xymatrix{
(x,\alpha)  \ar[rr]^{h'} \ar[d]_{\mathrm{Id}}& & W \ar[d]^{g}       \\
(x,\alpha)  \ar[r]_{h}& V \ar[r]_{f} &  X.
}$$
We have 
$$\begin{array}{ll}
(\Lambda^1 \otimes iso_{\alpha})(T) \circ \mu_V ([(x,\alpha) \xrightarrow{h} V])&=(\Lambda^1 \otimes
iso_{\alpha})(T)(h(x) \otimes [(x,\alpha) \xrightarrow{h} V])
\end{array}$$
$$
                                          =\left\lbrace
\begin{array}{cc}
p_g \circ f \circ h(x) \otimes [(x,\alpha) \xrightarrow{h'} W]&\mathrm{if\ } P =(x,\alpha)\\
0 & \mathrm{otherwise}
\end{array}
\right.$$
and
$$\begin{array}{ll}
\mu_W \circ iso_{\alpha}(T)( [(x,\alpha) \xrightarrow{h} V])&=\mu_W \left\lbrace
\begin{array}{cc}
 [(x,\alpha) \xrightarrow{h'} W]&\mathrm{if\ } P =(x,\alpha)\\
0 & \mathrm{otherwise}
\end{array}
\right.\\
\end{array}$$
  $$=\left\lbrace
\begin{array}{cc}
h'(x) \otimes [(x,\alpha) \xrightarrow{h'} W] &\mathrm{if\ } P =(x,\alpha)\\
0 & \mathrm{otherwise.}
\end{array}
\right.$$
By commutativity of the previous cartesian diagram, we have
$g \circ h'= f \circ h$ hence, by composition with $p_g$, we obtain: $h'=p_g \circ f \circ h$. Consequently, the linear maps $\mu_W$ give rise to a nonzero natural map $\mu: iso_{\alpha} \rightarrow \Lambda^1 \otimes iso_{\alpha}$. We deduce from the simplicity of the functors $iso_{\alpha}$ given in Corollary \ref{1.8} that the natural map $\mu$ is a monomorphism in \Fq.
\item
We check the commutativity of the following diagram, for a morphism $T=[V \xrightarrow{f} X \xleftarrow{g} W ]$ of $\mathrm{Hom}_{\Tq}(V,W)$:
$$\xymatrix{
(\Lambda^1 \otimes iso_{\alpha})(V)\ar[r]^-{\nu_V} \ar[d]_{(\Lambda \otimes iso_{\alpha})(T)}& iso_{\alpha}(V)   \ar[d]^{iso_{\alpha}(T)}       \\
 (\Lambda^1 \otimes iso_{\alpha})(W) \ar[r]_-{\nu_W} & iso_{\alpha}(W).
}$$
With the same notations as above, we have:
$$\begin{array}{ll}
 iso_{\alpha}(T) \circ \nu_V (v \otimes [h])&=iso_{\alpha}(T)(B(h(x),v)[h])\\
                                           &=\left\lbrace
\begin{array}{cc}
B(h(x),v) [h']&\mathrm{if\ } P \simeq\mathrm{Vect}( x )\\
0 & \mathrm{otherwise.}
\end{array}
\right.
\end{array}$$
and
$$\begin{array}{ll}
\nu_W \circ (\Lambda^1 \otimes  iso_{\alpha})(T)(v \otimes [h])&=\nu_W \left\lbrace
\begin{array}{cc}
p_g \circ f (v) \otimes [h']&\mathrm{if\ } P \simeq \mathrm{Vect}( x)\\
0 & \mathrm{otherwise.}
\end{array}
\right.\\
                         &=\left\lbrace
\begin{array}{cc}
B(h'(x), p_g \circ f (v))[h']&\mathrm{if\ } P \simeq\mathrm{Vect}( x )\\
0 & \mathrm{otherwise.}
\end{array}
\right.
\end{array}$$
We have:
$$\begin{array}{lll}
B(h(x),v) &= B(f \circ h(x), f(v)) &\mathrm{since\ } f
\mathrm{\ preserves\ quadratic\ forms}\\
          &= B(g \circ h'(x), f(v)) &\mathrm{by\
            commutativity\ of\ the\ cartesian\ diagram}\\
          &=B(g \circ h'(x), g \circ p_g \circ f(v))
          &\mathrm{by\ orthogonality\ of\ }W \mathrm{\ and\ }
          L\\
          &=B( h'(x), p_g \circ f(v))
          &\mathrm{since\ }g \mathrm{\ preserves\ quadratic\ forms.}\\
\end{array}$$
Hence, the linear maps $\nu_W$ define a natural map $\nu: \Lambda^1 \otimes iso_{\alpha} \rightarrow iso_{\alpha} $,  which is clearly nonzero. We deduce from the simplicity of the functors $iso_{\alpha}$ given in Corollary \ref{1.8}, that $\nu$ is an epimorphism of \Fq.
\end{enumerate}
\end{proof}
The natural maps $\mu$ and $\nu$ defined in the previous lemma allow us to define the following morphisms in \Fq.
\begin{defi} \label{6.2}
Let $n$ be a non-negative integer.
\begin{enumerate}
\item
The natural map $\mu_n: \Lambda^n \otimes iso_{\alpha} \rightarrow
\Lambda^{n+1} \otimes iso_{\alpha}$ is obtained by the following composition
$$\Lambda^n \otimes iso_{\alpha} \xrightarrow{1 \otimes \mu} \Lambda^n \otimes
\Lambda^1 \otimes iso_{\alpha} \xrightarrow{m \otimes 1} \Lambda^{n+1} \otimes
iso_{\alpha}$$
where $m:\Lambda^n \otimes \Lambda^1 \rightarrow \Lambda^{n+1}$ is the product in the exterior algebra.
\item
The natural map $\nu_n: \Lambda^{n+1} \otimes iso_{\alpha} \rightarrow
\Lambda^{n} \otimes iso_{\alpha}$ is obtained by the following composition 
$$\Lambda^{n+1} \otimes iso_{\alpha} \xrightarrow{\Delta \otimes 1} \Lambda^n \otimes
\Lambda^1 \otimes iso_{\alpha} \xrightarrow{1 \otimes \nu} \Lambda^{n} \otimes
iso_{\alpha}$$
where $\Delta:\Lambda^{n+1}  \rightarrow \Lambda^{n} \otimes \Lambda^1$ is the coproduct. 
%from est
%l'application coproduit de $\Lambda^{n+1}$ dans $\Lambda^{n}
%\otimes \Lambda^1$ dont on rappelle la définition: $$\Delta_V(v_1 \wedge \ldots \wedge v_{n+1})=\sum_i ((v_1
%\wedge \ldots \wedge \hat{v_i} \wedge \ldots \wedge v_{n+1}) \otimes
%v_i) $$
%pour un object $V$ de $\Tq$ et $v_1, \ldots, v_{n+1}$ des éléments de $V$.
\end{enumerate}
\end{defi}
We have the following proposition.
\begin{prop} \label{6.3}
The following sequence is an exact complex.
$$\ldots \rightarrow \Lambda^{n} \otimes iso_{\alpha} \xrightarrow{\mu_n}
\Lambda^{n+1} \otimes iso_{\alpha} \xrightarrow{\mu_{n+1}} \Lambda^{n+2}
\otimes iso_{\alpha} \rightarrow
\ldots$$
\end{prop}
\begin{proof}
We prove, according to the definition, that the kernel of $({\mu_{n+1}})_V$ is the vector space generated by the set
$$\{ v_1 \wedge \ldots \wedge v_n \wedge h(x) \otimes[h] \mathrm{\
  for\ } [h] \mathrm{\ a\ generator\ of\ }
  iso_{\alpha}(V)\mathrm{\ and\ } v_1, \ldots, v_n \mathrm{\
  elements\ in\ } V \}$$
and this space coincide with the image of $({\mu_{n}})_V$.
\end{proof}
Proposition \ref{6.3} justifies the introduction of the following functor.
\begin{defi} 
The functor $K_{\alpha}^n$ is the kernel of the map:
$\mu_n: \Lambda^{n}
\otimes iso_{\alpha} \rightarrow \Lambda^{n+1} \otimes iso_{\alpha}$.
\end{defi}
As observed in the proof of Proposition \ref{6.3}, we have the following characterization of the spaces $K_{\alpha}^n(V)$.
\begin{lm} \label{6.5}
For an object $V$ in $\Tq$, the space $K_{\alpha}^n(V)$ is generated by the following elements
$$z \wedge h(x) \otimes [h]$$
where $[h]$ is a canonical generator of the space $iso_{\alpha}(V)$ and $z$ is an element in $\Lambda^{n-1}(V).$
\end{lm}

The result below is a straightforward consequence of Proposition \ref{6.3}.
\begin{cor}
For $n$ a nonzero integer, we have the short exact sequence:
$$0 \rightarrow K_{\alpha}^n \rightarrow \Lambda^{n} \otimes iso_{\alpha} \rightarrow
K_{\alpha}^{n+1} \rightarrow 0.$$
\end{cor}
%-----------------------------------------------------------
\begin{rem}
We will prove, in a subsequent paper concerning the calculation of $\mathrm{Hom}_{\Fq}(\Lambda^n \otimes iso_{\alpha} , \Lambda^m \otimes iso_{\alpha} )$ that this short exact sequence is not split. \end{rem}
%-------------------------------------------------------------
%\aside{cette s.e.c. n'est pas scindée contrairement à ce que
%  j'énonçais dans la première version.\\
%Montrer qu'elle n'est pas scindée?}

We next explain how to decompose the functors $K_{\alpha}^n$. We begin by investigating the case of the functor $K_{\alpha}^1$.

\begin{lm} \label{6.8}
The functor $K_{\alpha}^1$ is equivalent to the functor $iso_{\alpha}$.
\end{lm}
\begin{proof}
Let $V$ be an object in $\Tq$. A basis of the vector space $K_{\alpha}^1(V)$ is given by the set of elements of the following form: $h(x) \otimes [h]$ for $[h]$ a canonical generator of $iso_{\alpha}(V)$. Then we define the linear map
$$\begin{array}{lll}
K_{\alpha}^1(V) & \xrightarrow{\sigma_V} & iso_{\alpha}(V) \\
h(x) \otimes [h] & \longmapsto & [h]
\end{array}$$
and we leave the reader to check that $\sigma_V$ is an isomorphism and that these linear maps are natural.
\end{proof}
In order to identify the composition factors of the functors $K_{\alpha}^n$ for $n>1$, we need the following lemma.
\begin{lm} \label{6.9}
For $n$ a nonzero integer, the morphism $\nu_n:\Lambda^{n+1} \otimes iso_{\alpha} \rightarrow
\Lambda^{n} \otimes iso_{\alpha}$ induces a morphism $\nu_n^K: K_{\alpha}^{n+1}
\rightarrow K_{\alpha}^n$ making the following diagram commute
$$\xymatrix{
K_{\alpha}^{n+1} \ar@{^{(}->}[d] \ar[r]^{\nu_n^K} & K_{\alpha}^n \ar@{^{(}->}[d]\\
\Lambda^{n+1} \otimes iso_{\alpha} \ar[r]_{\nu_n}& \Lambda^{n} \otimes iso_{\alpha}.
}$$
\end{lm}
\begin{proof}
For an object $V$ in $\Tq$ and $v_1 \wedge \ldots \wedge v_n
\wedge h(x) \otimes [h]$ a generator of $K_{\alpha}^{n+1}(V)$, we have:
$$({\nu_n}_V)(v_1 \wedge \ldots \wedge v_n
\wedge h(x) \otimes [h])$$
$$=(\sum_{i=1}^n v_1 \wedge \ldots \wedge \hat{v_i} \wedge \ldots  \wedge v_n
\wedge h(x) \otimes B(v_i, h(x)) [h]) + v_1 \wedge \ldots \wedge v_n
\otimes B(h(x),h(x))[h].$$
Since $B$ is alternating, we have $B(h(x),h(x))=0$. Hence,
$$({\nu_n}_V)(v_1 \wedge \ldots \wedge v_n
\wedge h(x) \otimes [h])\\
=\sum_{i=1}^n v_1 \wedge \ldots \wedge \hat{v_i} \wedge \ldots  \wedge v_n
\wedge h(x) \otimes B(v_i, h(x)) [h] \in K_{\alpha}^n(V).$$
We deduce the existence of the induced morphism $\nu_n^K$.
\end{proof}

This lemma justifies the introduction of the following definition.
\begin{defi} \label{6.10} 
For $n \ge 2$ an integer, let $L_{\alpha}^{n}$ be the kernel of the morphism $\nu_{n-1}^K: K_{\alpha}^{n}
\rightarrow K_{\alpha}^{n-1}$.
\end{defi}
We have the following characterization of the spaces $L_{\alpha}^n(V)$
which is useful below.
\begin{lm} \label{6.11}
For an object $V$ in $\Tq$, the space $L_{\alpha}^n(V)$ is generated by the elements of the following form:
$$z \wedge h(x) \otimes [h],$$
where $[h]$ is a canonical generator of the space $iso_{\alpha}(V)$ and $z$ is an element of $\Lambda^{n-1}(\mathrm{Vect}( h(x))^{\bot}).$
\end{lm}
\begin{proof}
Let $V$ be an object in $\Tq$. Since the space $L_{\alpha}^n(V)$ is a subvector space of $K_{\alpha}^n(V)$, we deduce from Lemma \ref{6.5} that the vector space  $L_{\alpha}^n(V)$ is generated by the elements of the following form: $z
\wedge h(x) \otimes [h]$.

For a given canonical generator $[h]$ of $iso_{\alpha}(V)$, since the quadratic space $V$
is nondegenerate, there is an element $w$ in $V$ such that  
$$B(h(x),w)=1.$$
Let $W$ be the space $\mathrm{Vect}(h(x),w)$ and $V \simeq W \bot W^{\bot}$ be an orthogonal decomposition of the space $V$. 
The canonical generator $z \wedge h(x) \otimes [h]$ of $K_{\alpha}^n(V)$ can be written in the following form:
$$z' \wedge h(x) \otimes [h]+z'' \wedge w \wedge h(x) \otimes [h] $$
where $z' \in \Lambda^{n-1}(W^{\bot})$ and $z'' \in \Lambda^{n-2}(W^{\bot})$.

Let $x$ be an element of $L_{\alpha}^n(V)$. Then:
\begin{equation*}
x=\sum_{[h] \in \mathcal{G}(iso_{\alpha}(V))} z_h \wedge h(x) \otimes [h]=\sum_{[h] \in \mathcal{G}(iso_{\alpha}(V))} (z'_h \wedge h(x) \otimes [h]+z''_h \wedge w \wedge h(x) \otimes [h])
\end{equation*}
where $\mathcal{G}(iso_{\alpha}(V))$ is the set of the canonical generators of $iso_{\alpha}(V)$. Consequently, 
$$\nu_{n-1}^K(z'_h \wedge h(x) \otimes [h])=0$$
since $z'_h \in \Lambda^{n-1}(W^{\bot}) \subset \Lambda^{n-1}(
\mathrm{Vect}( h(x))^{\bot}) $ and 
$$\nu_{n-1}^K(z''_h \wedge w \wedge h(x) \otimes [h])=z''_h \wedge
h(x) \otimes [h].$$
Since the element $x$ is in the kernel of $\nu_{n-1}^K$, we deduce that $z''_h=0$. Hence 
\begin{equation*}
x=\sum_{[h] \in \mathcal{G}(iso_{\alpha}(V))} z'_h \wedge h(x) \otimes [h]
\end{equation*}
for $z'_h \in \Lambda^{n-1}(W^{\bot}) \subset \Lambda^{n-1}(
\mathrm{Vect}( h(x) )^{\bot}) $.
\end{proof}
We have the following lemma.
\begin{lm} \label{6.12}
The composition: $ \Lambda^{n+2} \otimes iso_{\alpha} \xrightarrow{\nu_{n+1}} \Lambda^{n+1}
\otimes iso_{\alpha} \xrightarrow{\nu_{n}} \Lambda^{n} \otimes iso_{\alpha}$ is zero.
\end{lm}
\begin{proof}
Let $V$ be an object in \Tq\ and $v_1 \wedge \ldots \wedge v_{n+2} \otimes
[h]$ be an element of $\Lambda^{n+2} \otimes iso_{\alpha}(V)$. Then
$$\begin{array}{l}
\nu_n \circ \nu_{n+1} (v_1 \wedge \ldots \wedge v_{n+2} \otimes
[h]) \\
=\nu_n(\sum_{i=1}^{n+2}(v_1 \wedge \ldots \wedge \hat{v_i}
\wedge \ldots \wedge v_{n+2} \otimes B(v_i, h(x)) [h])\\
=\sum_{j \ne i} \sum_{i=1}^{n+2}(v_1 \wedge \ldots \wedge \hat{v_i}
\wedge \ldots \wedge \hat{v_j}
\wedge\ldots \wedge v_{n+2} \otimes B(v_i, h(x)) B(v_j,h(x)) [h])\\
=0 
\end{array}$$
since the characteristic is equal to $2$.
\end{proof}
We deduce the following result.
\begin{lm} \label{6.13}
The map $\nu_{n}^K: K_{\alpha}^{n+1} \rightarrow K_{\alpha}^{n}$ factors through $L_{\alpha}^n$.
\end{lm}
\begin{proof}
By Lemma \ref{6.9}, the following diagram is commutative:
$$\xymatrix{
K_{\alpha}^{n+1} \ar@{^{(}->}[d] \ar[r]^{\nu_{n}^K} &K_{\alpha}^{n} \ar@{^{(}->}[d] \ar[r]^{\nu_{n-1}^K} & K_{\alpha}^{n-1} \ar@{^{(}->}[d]\\
\Lambda^{n+1} \otimes iso_{\alpha} \ar[r]_{\nu_{n}}&\Lambda^{n} \otimes iso_{\alpha} \ar[r]_-{\nu_{n-1}}& \Lambda^{n-1} \otimes iso_{\alpha}.
}$$
Consequently, we deduce from Lemma \ref{6.12}, that $\nu_{n-1}^K \circ
\nu_{n}^K=0$. Hence, there is a morphism $\tilde{\nu}_n^K:K_{\alpha}^{n+1}
      \rightarrow L_{\alpha}^{n}$ making the following diagram commutative
      $$\xymatrix{
      &K_{\alpha}^{n+1} \ar@{.>}[ld]     \ar[d]_{\nu_{n}^K} \ar[dr]^{0}\\
L_{\alpha}^{n}=\mathrm{Ker}(\nu_{n-1}) \ar[r] & K_{\alpha}^{n} \ar[r]_{\nu^{K}_{n-1}} &K_{\alpha}^{n-1}.
}$$
\end{proof}
Then, we have the following proposition.
\begin{prop} \label{6.14}
For $n$ a nonzero integer, there is a short exact sequence:
$$0 \rightarrow L_{\alpha}^{n+1} \rightarrow K_{\alpha}^{n+1} \rightarrow
L_{\alpha}^{n} \rightarrow 0.$$
\end{prop}

\begin{proof}
It is sufficient to prove that the natural map $\tilde{\nu}_n^K:K_{\alpha}^{n+1} \rightarrow L_{\alpha}^{n}$ constructed in the proof of the previous Lemma is an epimorphism of \Fq.

Let $V$ be an object in \Tq and let $v_1 \wedge \ldots \wedge v_{n-1} \wedge
h(x) \otimes [h]$ be a generator of $L_{\alpha}^n(V)$. By definition of $L_{\alpha}^n(V)$
we have $B(v_i,h(x))=0$ for all $i$ in $\{1, \ldots, n-1 \}$. Since $h(x)$ is a nonzero element in the nondegenerate quadratic space $V$, there is an element $v$ in $V$ such that
$B(v,h(x))=1$. Then, we prove that the element  
$$v_1 \wedge \ldots \wedge v_{n-1} \wedge v \wedge h(x) \otimes [h]
\mathrm{\ de\ }K_{\alpha}^{n+1}(V)$$
verifies
$$(\tilde{\nu}_n^K)_V(v_1 \wedge \ldots \wedge v_{n-1} \wedge v \wedge h(x)
\otimes [h])=v_1 \wedge \ldots \wedge v_{n-1} \wedge h(x) \otimes
[h].$$
Hence, $\tilde{\nu}_n^K$ is surjective.
\end{proof}

\begin{rem}
Proposition \ref{6.14} is equivalent to the following statement: the complex 
$$ \ldots \rightarrow K_{\alpha}^{n+1} \xrightarrow{\nu_{n}^K} K_{\alpha}^{n} \xrightarrow{\nu_{n-1}^K} K_{\alpha}^{n-1} \rightarrow \ldots$$
is exact.
\end{rem}

\subsection{Simplicity of the functors $L_{\alpha}^n$}
In this section, we prove the following result, where the functors $L_{\alpha}^n$ are the subfunctors of  $\Lambda^n \otimes iso_{\alpha}$ defined in Definition \ref{6.10}.

\begin{thm} \label{6.15}
The functors $L_{\alpha}^n$ are simple.
\end{thm}

To prove this theorem, we need the following fundamental lemma.

\begin{lm} \label{6.16}
If $J$ is a subfunctor of $L_{\alpha}^n$, then for all object $V$ in
\Tq,  either $J(V)=\{0\}$, or $J(V)=L_{\alpha}^n(V)$.
\end{lm}

\begin{proof}
Let $J$ be a subfunctor of $L_{\alpha}^n$ and $V$ be an object in \Tq. Suppose that $J(V) \ne \{0\}$ and denote by $y$ a nonzero element of $J(V)$. We have 
\begin{equation} \label{Ln-E1}
y=\sum_{[h] \in \mathcal{G}(iso_{\alpha}(V))} z_h \wedge h(x) \otimes [h]
\end{equation}
where $z_h$ is an element of $\Lambda^{n-1}(\mathrm{Vect}( h(x))^{\bot})$ by Lemma \ref{6.11} and
$\mathcal{G}(iso_{\alpha}(V))$ is the set of canonical generators of the space $iso_{\alpha}(V)$.

The proof is divided into three steps; in the first one we prove that there exists a generator $[h]$ of $iso_{\alpha}(V)$ such that $z_h
\wedge h(x) \otimes [h] \in J(V)$. We deduce, in the second part, that for all other generator $[h']$ of $iso_{\alpha}(V)$ we have a nonzero element of the form $z'
\wedge h'(x) \otimes [h']$ in $J(V)$. Finally, we prove that for each element of the form $v_1 \wedge \ldots \wedge v_{n-1}$ in
$\Lambda^{n-1}( \mathrm{Vect}( h(x) )^{\bot})$ and each canonical generator $[h]$ of
$iso_{\alpha}(V)$, the element $v_1 \wedge \ldots \wedge v_{n-1} \wedge h(x)
\otimes [h]$ belongs to $J(V)$. This will prove that the two spaces $J(V)$ and $L_{\alpha}^n(V)$ are
isomorphic, by the characterization of the space $L_{\alpha}^n(V)$ given in Lemma \ref{6.11}.

\begin{enumerate}
\item
Let $[h]$ be a canonical generator of $iso_{\alpha}(V)$ such that, in the decomposition of $y$ given in  (\ref{Ln-E1}) the element $z_h \wedge h(x) \otimes [h]$ is nonzero. Since the space $V$ is nondegenerate, there exists an element $v$ in $V$ such that $B(v,h(x))=1$. We deduce a symplectic decomposition of $V$ of the following form:
\begin{equation} \label{Ln-E2}
V=\mathrm{Vect}(h(x),v) \bot \mathrm{Vect}(v_1,w_1) \bot \ldots \bot
\mathrm{Vect}(v_m,w_m)= \mathrm{Vect}(h(x),v) \bot V'.
\end{equation}
We consider the morphism of $\Eq$:
$$\begin{array}{llll}
f:& V& \rightarrow & V \bot (H_0)^{\bot (2m+1)}\\
  & h(x)& \longmapsto & h(x) \\
  & v& \longmapsto & v+a_0^1\\
%  &  &  \ldots   &      \\
  & v_k &\longmapsto & v_k +a_0^{2k}\\
 & w_k &\longmapsto & w_k +a_0^{2k+1}\\
% &  &  \ldots   &      \\
\end{array}$$
for $k$ an integer between $1$ and $m$, which allows us to define the following morphism in $\Tq$: 
$\xymatrix{
T=[V \ar[r]^(.4){f} & V \bot (H_0)^{\bot (2m+1)} & \ar@{_{(}->}[l]_(.3){i} V]}.$ 
We deduce from the two cartesian diagrams below:
$$\xymatrix{
(x,\alpha) \ar[d]_{\mathrm{Id}} \ar[rr]& & V \ar@{^{(}->}[d]^i\\
(x,\alpha) \ar[r]_{h}& V \ar[r]_(.3){f} &  V \bot (H_0)^{\bot (2m+1)}}$$
 and, for $h_i \ne h$,
$$\xymatrix{
\{ 0 \} \ar[d]_{\mathrm{Id}} \ar[rr]& & V \ar@{^{(}->}[d]^i\\
(x,\alpha) \ar[r]_{h_i}& V \ar[r]_(.3){f} &  V \bot (H_0)^{\bot (2m+1)}}$$
that  $J(T)(y)=z_h \wedge h(x) \otimes [h] \in J(V)$,
since  $\epsilon(T)=\mathrm{Id}_V$.
\item
Let $[h']$ be a canonical generator of $iso_{\alpha}(V)$ different from $[h]$. We have
$q(h'(x))=q((h(x))=\alpha$, hence the linear isomorphism denoted by $\underline{f}$, from $(h(x),\alpha)$ to $(h'(x),\alpha)$ is a morphism in $\Eqd$. So, we can apply Theorem \ref{1.3}
to obtain the existence of a morphism $f$ of $\mathrm{Hom}_{\Eq}(V, V)$ making the following diagram commutative:
$$\xymatrix{
V  \ar[r]^f & V  \\
(h(x), \alpha) \ar@{^{(}->}[u] \ar[r]_{\underline{f}}& (h'(x), \alpha). \ar@{^{(}->}[u] 
}$$
We deduce the following cartesian diagram: 
$$\xymatrix{
(x,\alpha) \ar[d]_{\mathrm{Id}} \ar[rr]^{h'}& & V \ar[d]^{\mathrm{Id}}\\
(x,\alpha) \ar[r]_{h}& V \ar[r]_{f} &  V. }$$
Consequently, by the consideration of the morphism $T=[V \xrightarrow{f} V
\xleftarrow{\mathrm{Id}} V]$, we obtain:
$$J(T)(z_h \wedge h(x) \otimes [h])=\Lambda^{n-1} (f)(z_h) \wedge h'(x)
\otimes [h'] \in J(V).$$

\item
By the point ($1$) of the proof, there exists a nonzero element of the form $z \wedge h(x) \otimes [h]$ in $J(V)$. We want to prove that, for each element of the form $v_1 \wedge \ldots \wedge v_{n-1}$ in
$\Lambda^{n-1}(\mathrm{Vect}( h(x) )^{\bot})$ and each canonical generator $[h]$ of
$iso_{\alpha}(V)$, the element $v_1 \wedge \ldots \wedge v_{n-1} \wedge h(x)
\otimes [h]$ belongs to $J(V)$. According to the proof of Lemma \ref{6.11}, it is sufficient to prove that the element $v_1 \wedge \ldots \wedge v_{n-1} \wedge h(x)
\otimes [h]$ belongs to $J(V)$ for $v_1 \wedge \ldots \wedge v_{n-1}$
in $\Lambda^{n-1}(V')$ where $V'$ is the space considered in the decomposition (\ref{Ln-E2}). By simplicity of the functor $\Lambda^{n-1}$ in $\F$, we have the existence of an endomorphism $g$ of $\epsilon(V')$ such that 
$$\Lambda^{n-1} (g)(z)=v_1
\wedge \ldots \wedge v_{n-1}.$$
We deduce that
$$J(\mathrm{Id}_{\Tq}(\mathrm{Vect}(h(x),v)) \bot t_{g})(z \wedge h(x) \otimes [h]) = v_1
\wedge \ldots \wedge v_{n-1} \wedge h(x)
\otimes [h]$$
where $t_{g}$ is an antecedent of $g \in \mathrm{End}_{\E^f}(\epsilon(V'))$ by the forgetful functor $\epsilon: \Tq \rightarrow \E^f$ which is full by Proposition $3.5$ in \cite{math.AT/0606484} and $\mathrm{Id}_{\Tq}(\mathrm{Vect}(h(x),v))
\bot t_{g}$ is the orthogonal sum of the morphisms of \Tq.
\end{enumerate}
\end{proof}
\begin{proof}[Proof of Theorem \ref{6.15}]
Let $J$ be a nonzero subfunctor of $L_{\alpha}^n$ and $V$ be an object in $\Tq$
such that the space $J(V)$ is nonzero. By Lemma \ref{6.16}, we have $J(V)=L_{\alpha}^n(V)$. We prove, in the following, that for all object $W$ in $\Tq$, $J(W)=L_{\alpha}^n(W)$.

Let $W$ be a fixed object in $\Tq$. The proof is divided into two parts.
\begin{enumerate}
\item{Let us prove that $J(V \bot W) \simeq L_{\alpha}^n(V \bot W)$.}

Let $i$ be the canonical inclusion from $V$ to $V \bot W$ and 
$T=[V \xrightarrow{i} V \bot W \xleftarrow{\mathrm{Id}} V \bot W]$ be the morphism of $\Tq$. We have the following cartesian diagram:
$$\xymatrix{
(x,\alpha) \ar[d]_{\mathrm{Id}} \ar[rr]^{i \circ h}& & V \bot W \ar[d]^{\mathrm{Id}}\\
(x,\alpha) \ar[r]_{h}& V \ar[r]_{i} & V \bot  W.}$$
Since 
$$[V \bot W \xrightarrow{\mathrm{Id}} V \bot W \xleftarrow{i} V ] \circ [V \xrightarrow{i} V \bot W \xleftarrow{\mathrm{Id}} V \bot W]=Id_V$$
we deduce that the space $J(V \bot W)$ is nonzero. Hence, by Lemma \ref{6.16}, we have $J(V \bot W) \simeq L_{\alpha}^n(V \bot W)$.
\item{Let us prove that $J(W) \simeq L_{\alpha}^n(W)$.}

According to Lemma \ref{6.16}, it is sufficient to prove that if $L_{\alpha}^n(W)$ is 
non trivial then $J(W)$ is not zero.
Let $j$  be the canonical inclusion from $W$ to $V \bot W$.

We deduce from 
$$[V \bot W \xrightarrow{\mathrm{Id}} V \bot W \xleftarrow{j} W ] \circ [W \xrightarrow{j} V \bot W \xleftarrow{\mathrm{Id}} V \bot W]=Id_W$$
the existence of the surjection:
$\xymatrix{
L_{\alpha}^n(V \bot W) \ar@{->>}[r] & L_{\alpha}^n(W)}.$
Furtermore, by the first point of the proof, we have $J(V \bot
W)=L_{\alpha}^n(V \bot W)$. This gives rise to the following commutative diagram:
$$\xymatrix{
J(V \bot W) \ar[r] \ar[d]_{\simeq} & J(W) \ar@{^{(}->}[d] \\
L_{\alpha}^n(V \bot W) \ar@{->>}[r] & L_{\alpha}^n(W).}$$
Consequently, by a diagram chasing argument, we obtain that if $L_{\alpha}^n(W)$ is nonzero then
$J(W)$ is nonzero.
\end{enumerate}

\end{proof}

We prove in the following proposition that this gives rise to two families of non isomorphic functors. \begin{prop}
The functors in the union of the families $\{L^n_0\  |\  n \in \mathbb{N} \}$ and
$\{L^n_1\ |\ n \in \mathbb{N}\}$ are pairwise non-isomorphic.
\end{prop}

\begin{proof}
For a fixed $\alpha$, there exists, for all integer $n$, a minimal integer $d(n)$ such that  
$$L^n_{\alpha}(H_0^{\bot d(n)}) \ne 0.$$
Let $k$ be an integer different from $n$, if $|n-k| \geq 2$, the integers $d(n)$ and $d(k)$ allow us to distinguish the simple functors $L^n_{\alpha}$ and $L^k_{\alpha}$, in the contrary case, we prove that the dimensions of the spaces $L^n_{\alpha}(H_0^{\bot d(n)})$ and
$L^k_{\alpha}(H_0^{\bot d(n)})$  are different. This proves that the simple functors $L^n_{\alpha}$ and $L^k_{\alpha}$ are not isomorphic.

Furthermore, two simple functors $S_1$ and $S_2$ in $\Fq$ are not isomorphic if there exists a morphism $T$ in $\Tq$ such that $S_1(T)=0$ and $S_2(T) \ne 0$. Moreover, the morphisms $T$ constructed in the first point of the proof of Lemma \ref{6.16} verify 
$$L^n_{\alpha}(T) \ne 0 \mathrm{\quad and\quad } L^k_{(\alpha +1)}(T) = 0$$
where $(\alpha +1)$ is the reduction mod $2$ of $\alpha +1$.

\end{proof}

%\subsubsection{Quelques exemples de décompositions des foncteurs   $\Lambda^n \otimes iso_{\alpha}$}
%%%%%%%%%%%%%%%%%%%%%%%%%%%%%%%%%%%%%%%%%%%%%%%%%%%
\section{The composition factors of the functors $\M{0,1}$ and  $\M{1,1}$}
We prove in this section that the functors $\M{0,1}$ and  $\M{1,1}$ are uniserial (i.e. the lattice of its subfunctors is totally ordered) and the composition factors of the functor $\M{0,1}$ are the functors $L^n_0$ and those of $\M{1,1}$ are the functors $L^n_1$. For that, we identify the subquotient  $k_d \m{\alpha,1}/k_{d+1} \m{\alpha,1}$ of the filtration of the functor $\m{\alpha,1}$ by the functors $k_d \m{\alpha,1}$ introduced in Proposition \ref{5.2}
with the functor $L^{d+1}_{\alpha}$ defined in Definition \ref{6.10}. This gives rise to the following result.
\begin{thm} \label{7.1}
 The functor $ \m{\alpha,1}$ is uniserial and its unique composition series is given by the decreasing filtration by the functors $k_d \m{\alpha,1}$
$$\ldots \subset k_d \m{\alpha,1} \subset \ldots \subset k_1
\m{\alpha,1} \subset k_0 \m{\alpha,1}= \m{\alpha,1}$$
which verifies:
$$k_d \m{\alpha,1}/k_{d+1} \m{\alpha,1} \simeq L^{d+1}_\alpha.$$
\end{thm}

\begin{rem}
Remark that the strategy of the following proof of the uniseriality of $\m{\alpha,1}$ is close to the proof of Lemma \ref{6.16}.
\end{rem}
\begin{proof}
To prove that the functor $\m{\alpha,1}$ is uniserial, it is sufficient to prove that if $J$ is a nonzero subfunctor of  $\m{\alpha,1}$ then there exists an integer $d$ such that $J=k_d \m{\alpha,1}$.

Let $V$ be an object in $\Tq$ such that $J(V) \neq 0$ and $v$ be a nonzero element of $J(V)$, we have:
$$v=\sum_{ \{x,y\} \in A_V}  \alpha_{\{ x,y \} }(v)[ \{ x,y\} ]$$
where $A_V=\{ \{x,y\} | x \in V, y \in V, q(x+y)=\alpha, B(x,y)=1 \} $. One verifies easily that for $\{x,y\} \in A_V$ and $l \in V$, $\{x+l,y+l \} \in A_V$ if and only if $l \in (\mathrm{Vect}(x,y))^{\bot}$ or $l=x+y+l'$ where $l' \in (\mathrm{Vect}(x,y))^{\bot}$. 
Since $\{x+(x+y+l'), y+(x+y+l') \}=\{ y+l',x+l' \}=\{ x+l', y+l' \}$ we obtain, after reordering:
\begin{equation} \label{eq.7.1}
v=\sum_{i=1}^{n} \sum_{ l \in \mathcal{L}_i}  [ \{ x_i +l,y_i +l\} ]
\end{equation}
where $\forall i$ $\{x_i, y_i\} \in A_V$, for $i \neq j$, $x_i +y_i \neq x_j+y_j$ and  $\mathcal{L}_i$ is a subvector space of $(\mathrm{Vect}(x_i,y_i))^{\bot}$ of dimension $r_i$.  We consider a vector space $V'$, an element $v'$ of $J(V')$ and an element $\{ x', y'\}$ in $A_{V'}$ such that the dimension $r'$ of $\mathcal{L}'$ is minimal. We deduce from the decomposition \ref{eq.7.1} that $J \subset k_{r'} \m{\alpha,1}$. 

To prove that $k_{r'} \m{\alpha,1} \subset J$, we consider the following decomposition of $v' \in V'$:
$$v'=\sum_{l' \in \mathcal{L}'}[\{ x'+l', y'+l' \} ]+\sum_{j=1}^{p} \left( \sum_{ l \in \mathcal{L}_j} [ \{ x_j+l,y_j+l \}] \right)$$
where $\forall j\ x_j+y_j \neq x'+y'$.
By definition of $\m{\alpha,1}$ and $\M{\alpha,1}$, we have a natural map $\sigma: J \rightarrow P_{\mathbb{F}} \otimes iso_{\alpha}$ such that 
$$\sigma_{V'}(v')= \sum_{l' \in \mathcal{L}'} \left( [ x'+l'] +[y'+l' ]   \right) \otimes  [ x'+y']  +  \sum_{j=1}^{p}  \sum_{ l \in \mathcal{L}_j}  \left( [ x_j+l] +[y_j+l] \right) \otimes [x_j + y_j ]$$
where, by abuse, we denote by $v$ the linear map $\FF \rightarrow V$ determined by $v$.
Let $f: \epsilon(V') \rightarrow \epsilon(V')$ be the linear map such that $f(x'+y')=x'+y'$ and $f(x_j+y_j)=x_j+y_j+m$ where $m$ is a nonzero element of $\mathcal{L}_j$. Since $\epsilon$ is full by Propositon $3.5$ in \cite{math.AT/0606484}, we obtain the existence of a morphism $T$ in $\mathrm{Hom}_{\Tq}(V',V')$ such that $\epsilon(T)=f$.
We deduce that 
$$(P_{\mathbb{F}} \otimes iso_{\alpha})(T) \sigma_V(v')=\sum_{l' \in \mathcal{L}'} \left( [ x'+l'] +[y'+l' ]   \right) \otimes  [ x'+y'] $$
and, consequently
$$J(T)(v')=\sum_{l' \in \mathcal{L}'}[\{ x'+l', y'+l' \} ] \in J(V')$$
where $\mathcal{L}'$ is a subvector space of $(\mathrm{Vect}(x',y'))^{\bot}$ of dimension $r'$.

Let $W$ be an object of $\Tq$, $\sum_{l \in \mathcal{L}}[\{ w_1+l, w_2+l \} ]$, where $\mathcal{L}$ is a subvector space of $\mathrm{Vect}(w_1,w_2)^{\bot}$ of dimension $r'$,  be a generator of $k_{r'}\m{\alpha,1}(W)$ and $g:\epsilon(V') \rightarrow \epsilon(W)$ be the linear map such that $g(x'+y')=w_1+w_2$ and $g$ send a basis of $\mathcal{L}'$ to a basis of $\mathcal{L}$. By the fullness of $\epsilon$ we obtain a morphism $T'$ in $\mathrm{Hom}_{\Tq}(V',W)$ such that:
$$J(T')(\sum_{l' \in \mathcal{L}'}[\{ x'+l', y'+l' \} ])=\sum_{l \in \mathcal{L}}[\{ w_1+l, w_2+l \} ] \in J(W).$$
Hence $J=k_{r'}\m{\alpha,1}$.

By Lemma \ref{5.7}, we have:
$$(g_d \otimes iso_{\alpha}) \circ i_d(k_d \m{\alpha,1}) \subset
L_{\alpha}^{d+1}.$$
Consequently, by the commutative diagram (\ref{filt-diagramme}) given in the proof of the second point of Proposition \ref{5.5} we have the natural map
$$\sigma: k_d \m{\alpha,1}/k_{d+1} \m{\alpha,1} \rightarrow L_{\alpha}^{d+1}.$$
Since the quotients $k_d \m{\alpha,1}/k_{d+1} \m{\alpha,1}$ are nonzero by Proposition \ref{5.8} and the functors $L_{\alpha}^{d+1}$ are simple by Theorem \ref{6.15}, the natural map $\sigma$ is an equivalence.
\end{proof}

\bibliographystyle{amsplain}
\bibliography{these}

\end{document}